\documentclass[12pt]{amsart}
\usepackage{amssymb}
\usepackage{latexsym}
\usepackage{pstricks}
\usepackage{epsf}
\newtheorem{thm}{Theorem}[section] 
\newtheorem{pro}[thm]{Proposition}
\newtheorem{lem}[thm]{Lemma}
\newtheorem{cor}[thm]{Corollary}

\theoremstyle{definition}

\newtheorem*{note}{Note}

\theoremstyle{remark}
\newtheorem{rem}[thm]{Remark}

\def\1{{\rm1\mathchoice{\kern-0.25em}{\kern-0.25em}
        {\kern-0.2em}{\kern-0.2em}I}}

\newcommand{\one}{{ \rm \setlength{\unitlength}{1em}
\begin{picture}(0.7,0.65)
\put(0,0){$1$}\put(0.34,0){\line(0,1){0.65}}
\end{picture} }}
\newcommand{\lmn}[1]{\vadjust{\setbox1=\vtop{\hsize 25mm
\parindent=0pt\baselineskip=9pt
\rightskip=4mm plus 4mm#1}
\hbox{\kern-26mm\smash{\raise .5ex\box1}}}}

\newcommand{\psdiag}[3]{\hspace{1mm}\raisebox{-#1mm}{\epsfysize#2mm
\epsffile{#3.eps}}\hspace{1mm}}

\newcommand{\nc}{\newcommand}
\def\be#1\ee{\begin{equation}#1\end{equation}}
\def\bee#1\eee{\begin{equation*}#1\end{equation*}}
\nc{\bc}{\begin{center}}
\nc{\ec}{\end{center}}
\nc{\bb}{\mathbb}
\nc{\cal}{\mathcal}
\nc{\N}{\mathbb{N}}\nc{\Z}{\mathbb{Z}}
\nc{\oneun}{\one\!{}_1}

\def\Q{\mathbb{ Q}}
\def\K{{\cal K}}

\def\v8{\vskip 8pt}
\def\a{\alpha}

\def\la{\langle}
\def\ra{\rangle}
\def\l{\lambda}
\def\n{\nu}
\def\m{\mu}
\def\d{\delta}

\begin{document}

\title[Idempotents in BMW algebras]{Skein construction of idempotents in
Birman-Murakami-Wenzl algebras}
\author{Anna Beliakova}
\address{Mathematisches Institut, Rheinsprung 21,
CH-4051 Basel, Switzerland}
\email{Anna.Beliakova@unibas.ch}
\author{Christian Blanchet}
\thanks{The second author wishes to acknowledge the hospitality of
the Mathematisches Institut,
Universit\"at Basel.}
\address{LMAM, Universit\'e de Bretagne Sud, 1 rue de la Loi,
F-56000 Vannes, France}
\email{Christian.Blanchet$\char'100$univ-ubs.fr}

\begin{abstract}
We give skein theoretic formulas for minimal idempotents
in the Birman-Murakami-Wenzl algebras. These formulas
are then applied to derive various known results needed
in the construction of quantum invariants and modular categories.
 In particular, an elementary proof 
of the Wenzl formula  for quantum dimensions is given. 
This proof does not use the representation theory of quantum
groups and the character formulas. 
\end{abstract}

\maketitle
\section*{Introduction}
The Birman-Murakami-Wenzl algebras are deformations of the Brauer
centralizer algebras \cite{bw,Mu}. They are quotients of the
Artin braid groups algebras, and have appeared in connection with
the Kauffman link invariant and the quantum groups of types
B, C and D.
 The Birman-Murakami-Wenzl algebras are generically semi-simple, and their
structure was given by Wenzl \cite{Wbcd}.
 They play a key role in the construction of quantum invariants,
modular categories and Topological Quantum Field Theories, as was
shown by Turaev and Wenzl \cite{TW,TW2}.
 Our purpose here is to study the structure of these
algebras without using their representation theory.
 In a separate article, we will pursue Turaev and Wenzl's program
and construct four series of modular categories.
 Together with the present paper this construction
will be reasonably self-contained.

Our main results are the following.
\begin{itemize}
\item We give explicit formulas for minimal idempotents
in the Birman-Murakami-Wenzl algebras $K_n$. These are then used to obtain
the semi-simple decomposition of $K_n$, together with
a basis of {\it matrix units}. 
 Similar results were obtained by Ram and Wenzl  \cite{RW}
using Jones {\em basic construction}.

\item We give a skein theoretic proof of the Wenzl formula
for the quantum dimensions  of these minimal idempotents.
The key point is here the proof of the recursive formula (\ref{rec}).
This formula  is further used to  derive 
versions of the Wenzl formula corresponding to the quantum group
specializations and to discuss the  existence of  
idempotents in the non generic case.

\end{itemize} 

{\bf Conventions.}  Throughout this paper, the manifolds are 
 compact, smooth and oriented. 
By a {\it  link}
we mean an { isotopy class of an unoriented framed link}.
 Here, a framing is a non-singular normal vector field up
to homotopy. By a {\it  tangle} in a $3$-manifold $M$
we mean an { isotopy class of a framed tangle relative to the 
 boundary}.
Here the boundary of the tangle is a finite set
of points in $\partial M$, together with a non zero vector tangent to 
 $\partial M$ at each point.
 Note that a framing together with an orientation is equivalent
to a trivialization of the normal bundle up to homotopy.
By an {\it  oriented  link} 
we mean an isotopy class of a link  together with a 
trivialization of the normal bundle up to homotopy.
 By an {\it  oriented tangle}
we mean an { isotopy class of a tangle together with a 
trivialization of the normal bundle, up to homotopy
relative to the  boundary}.
Here the boundary of the tangle is a finite set
of points in $\partial M$,  together with a trivialization
of the tangent space to 
 $\partial M$ at each point.
 In the figures, a preferred convention using the plane
gives the framing (blackboard framing).

\section{The Birman-Murakami-Wenzl ribbon category}

\subsection{Kauffman skein relations.}
\begin{figure}[h]
$$
\psdiag{3}{9}{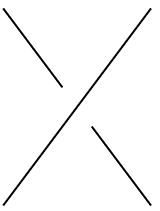}\;-\;\psdiag{3}{9}{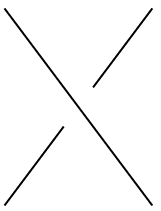} = \;(s-s^{-1})\;
\left(\;\,\psdiag{3}{9}{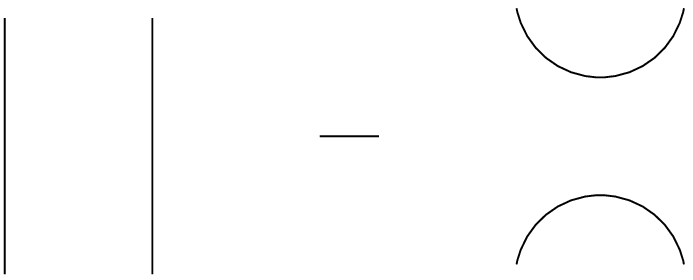}\;\,\right) 
$$
$$
\psdiag{3}{9}{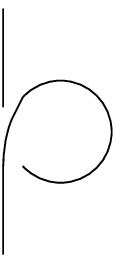}\;=\;\a\;\; \psdiag{3}{9}{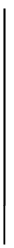}\;,\;\;\;\;\;
\psdiag{3}{9}{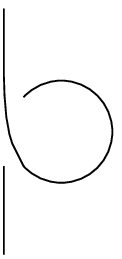}\;=\;\a^{-1} \;\;\psdiag{3}{9}{1}$$

$$L\;\amalg \;\bigcirc\; = \;\left(\frac{\a-\a^{-1}}{s-s^{-1}}+1
\right)\;\, L  $$
\caption{\label{kaufrel} Kauffman skein relations}
\end{figure}

Let  $M$ be a 3-manifold (possibly with a given finite set $l$ of  
framed points  on the boundary).
 We denote by $\K(M)$  (resp. $\K(M,l)$)
the $k$-module freely generated by
links in $M$ (and 
tangles in $M$ that meet $\partial M$ in $l$)
modulo (the relative isotopy  and)
the Kauffman  skein relations in
Figure~\ref{kaufrel}.

We suppose that $k$ is an integral domain 
containing the invertible elements $\a$, $s$ and that 
$\frac{\a-\a^{-1}}{s-s^{-1}}$ lives in $k$. 
We call $\K(M)$ the skein module  of $M$.

For example, $\K(S^3)\cong {k}$.
The isomorphism sends any link $L$
in $S^3$ to its
Kauffman polynomial $\langle L\rangle$, normalized by $\la\emptyset\ra=1$.

\noindent

\subsection{The Birman-Murakami-Wenzl  category.}
The Birman-Murakami-Wenzl
 category $K$ is
defined as follows. An object of $K$  is a disc
$D^2$ equipped with a finite set of points and a non zero tangent vector
at each point. Unless otherwise specified, we will use
 the second vector
of the standard basis (the vector $\sqrt{-1}$ in complex notation).
If $\beta=(D^2,l_0)$ and $\gamma=(D^2,l_1)$ are two objects,
the module $Hom_K(\beta,\gamma)$ is $\K(D^2\times [0,1],l_0\times 0\amalg
l_1\times 1)$. The notation $K(\beta,\gamma)$ and $K_\beta$ will be used for 
$Hom_K(\beta,\gamma)$ and $End_K(\beta)$, respectively.
For composition, we use the {\it covariant} notation:
$$
\begin{array}{rll}
K(\beta,\gamma)\times K(\gamma,\d)&\to& K(\beta,\delta)\\
 (f,g)&\mapsto &fg
\end{array}
$$
In our figures the time parameter goes upwards, so that the morphism $fg$
is depicted with $g$ lying above $f$, and the normal vector field
is orthogonal to the plane and points `inside the blackboard'. 

The Birman-Murakami-Wenzl
 category is a ribbon Ab-category (see \cite[Ch II]{Tu}).
 Ribbon categories admit a theory of traces of morphisms
and dimensions of objects for which we will use
the terminology {\it quantum trace} and {\it quantum dimension}.
 In the case of the category
of finite dimensional vector spaces, equipped with
trivial braiding and twist, these traces and dimensions coincide
with the usual ones
 \cite[Section I.1.7 and Lemma II.4.3.1]{Tu}.
 We will use the notation $\la f\ra $ for the
 quantum trace  of $f\in K_\beta$. This quantum
trace  is equal
to the value of the closure of $f$ in $ \K(S^3)\cong k$
obtained by gluing of $D^2\times \{0\}$
and  $D^2\times \{1\}$ in $K_\beta$ along the identity map.

We denote by $n$ the object formed by the $n$ points $\{ (2j-1)/n\; -1 \; :\;
j=1,...,n \}$ ($0$ is the trivial object).
 If we consider only these standard objects $n$, $n\geq 0$, 
we obtain a full subcategory equivalent to $K$, 
which was named the Kauffman category
by Turaev who first introduced it in \cite[Section 7.7]{T1}.
 
The  algebra $K_n=End_K(n)$ is
  isomorphic to the Birman-Murakami-Wenzl algebra 
which is
 the quotient of the braid group algebra
$k[B_n]$  by the Kauffman skein relations
\cite{BW,Mu}.   
 For a proof of the above isomorphism, see \cite{MT} or \cite{T1}.
This algebra
 is a deformation of the Brauer algebra (i.e. the centralizer
algebra of the semi-simple  Lie algebras of types B,C and D).
If $k$ is a field, then the  algebra $K_n$ is known to be 
semi-simple \cite{Wbcd}, except possibly if $s$ is a root of unity,
or $\alpha=\pm s^n$ for some $n\in\Z$. 
Its simple components correspond to the partitions $\l=(\l_1,
...,\l_p)$ with $|\l|=\sum_i \l_i=n-2r$, $r=0,1,...,[n/2]$.

\begin{figure}[h]
\begin{center}
\mbox{\epsfysize=3,5cm \epsffile{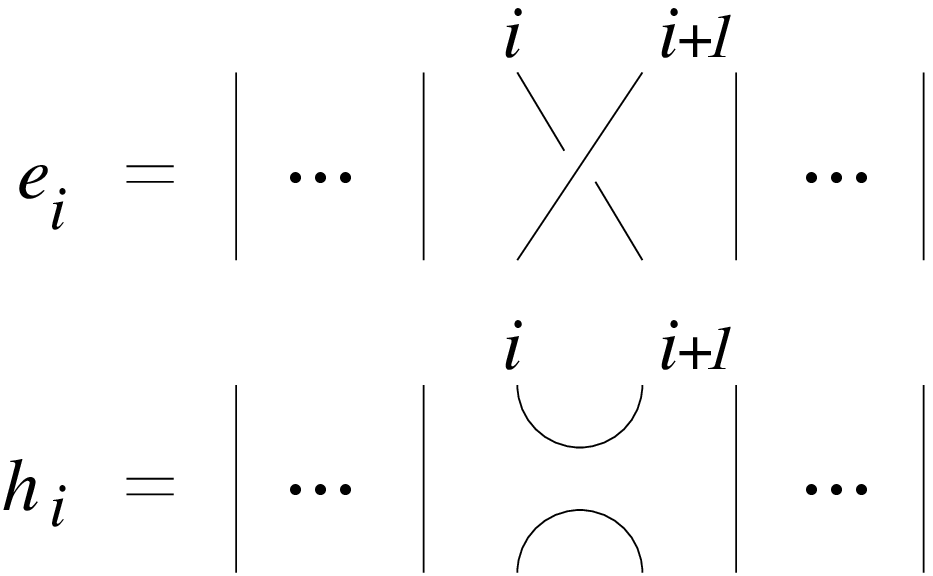}}
\end{center}
\caption{\label{generators} Generators of $K_n$}
\end{figure}

The Birman-Murakami-Wenzl 
algebra $K_n$ is generated by the identity
$\one\!{}_n$, positive transpositions  $e_1, ..., e_{n-1}$
and hooks $h_1,...,h_{n-1}$
 drawn in Figure~\ref{generators}.

If the $e_i$ are supposed to be invertible,
then a complete system of relations \cite{Wbcd} is given by
$$\begin{array}{crcl}
(B1)\ &e_ie_{i+1}e_i&=&e_ie_{i+1}e_i\\
(B2)\ &e_ie_j&=&e_je_i,\text{ for }|i-j|\geq 2\\
(R1)\ &h_ie_i&=&\a^{-1}h_i\\
(R2)\ &h_ie_{i-1}^{\pm 1}h_i&=&\a^{\pm 1}h_i\\
(K) \ &e_i-e_i^{-1}&=&(s-s^{-1})(\one\!{}_n-h_i)\ .
\end{array}$$
The quotient of $K_n$ by the ideal $I_n$ generated by $h_{n-1}$
is isomorphic to Hecke algebra ${H_n}$.
 We will use the knowledge of this
Hecke algebra $H_n$ to study $K_n$.
Note that $I_n=\{(a\otimes \one\!{}_1) h_{n-1} (b\otimes\one\!{}_1):
\; a,b \in K_{n-1}\}$.

\section{Hecke algebras} \label{hecke}
The  Hecke category $H$ is defined similarly as above,
using the Homfly skein theory.
An object in this category is a disc ${D}^2$ equipped with a set 
of points with a trivialization of the tangent space
at each point.
 If $\beta=({D}^2,l_0)$ and $\gamma=({D}^2,l_1)$
are two objects,
 the module $Hom_{{H}}(\beta,\gamma)={H}(\beta,\gamma)$ is
the Homfly skein module
$\mathcal{H}( {D}^2\times[0;1],l_0\times 0\amalg l_1\times 1)$.
 Here the Homfly skein module in $M$ is
freely generated by
oriented framed  tangles in $M$
modulo (the relative isotopy  and)
 the Homfly relations 
given in Figure~\ref{homflyrel}.
 Note that there
 we have specialized the three  variable  Homfly skein theory for 
framed links.
\begin{figure}[h]
\begin{center}
\vspace{5pt}
\begin{pspicture}[.4](0,0)(1,1)
\psline[linewidth=.5pt]{->}(.1,0)(.9,1)
\psline[linewidth=.5pt]{-}(.9,0)(.6,.4)
\psline[linewidth=.5pt]{->}(.4,.6)(.1,1)
\end{pspicture}
$\ -\ $
\begin{pspicture}[.4](0,0)(1,1)
\psline[linewidth=.5pt]{->}(.9,0)(.1,1)
\psline[linewidth=.5pt]{-}(.1,0)(.4,.4)
\psline[linewidth=.5pt]{->}(.6,.6)(.9,1)
\end{pspicture}
$\ =\ (s-s^{-1})$
\begin{pspicture}[.4](0,0)(1,1)
\pscurve[linewidth=.5pt]{->}(.1,0)(.4,.5)(.1,1)
\pscurve[linewidth=.5pt]{->}(.9,0)(.6,.5)(.9,1)
\end{pspicture}\\[10pt]
\begin{pspicture}[.4](0,0)(1,1)
\pscurve[linewidth=.5pt]{-}(.1,0)(.15,.45)(.45,.75)(.7,.5)(.4,.2)(.20,.38)
\pscurve[linewidth=.5pt]{->}(.14,.62)(.1,.9)(.1,1)
\end{pspicture}
$\ =\ \a$
\begin{pspicture}[.4](0,0)(1,1)
\psline[linewidth=.5pt]{->}(.5,0)(.5,1)
\end{pspicture}
\hspace{2cm}
\begin{pspicture}[.4](0,0)(1,1)
\pscurve[linewidth=.5pt]{<-}(.1,1)(.15,.55)(.45,.25)(.7,.5)(.4,.8)(.20,.62)
\pscurve[linewidth=.5pt]{-}(.14,.38)(.1,.1)(.1,0)
\end{pspicture}
$\ =\ \a^{-1}$
\begin{pspicture}[.4](0,0)(1,1)
\psline[linewidth=.5pt]{->}(.5,0)(.5,1)
\end{pspicture}\\[10pt]
$L\ \cup\ $
\begin{pspicture}[.4](0,0)(1,1)
\pscircle(.5,.5){.4}
\end{pspicture}
$\ =\ {\frac{\a-\a^{-1}}{s-s^{-1}}}\ L$
\vspace{5pt}
\end{center}
\caption{\label{homflyrel} Homfly skein relations}
\end{figure}

 We also simply denote  by $n$
the object formed by the $n$  
 points $\{ (2j-1)/n\; -1 \; :\;
j=1,...,n \}$,
 equipped with the standard trivialization.

The positive permutation braids $w_\pi$ represent
a basis of the module $H_n$, 
indexed by permutations $\pi$.
The symmetrizers  and antisymmetrizers
 in $H_n$  are represented respectively by the following
elements $f_n$ and $g_n$ of the braid group algebra.
$$f_n=\frac{1}{[n]!}s^{-\frac{n(n-1)}{2}}
\sum_{\pi\in \mathcal{S}_n} s^{l(\pi)}w_\pi\ ,$$
$$g_n=\frac{1}{[n]!}s^{\frac{n(n-1)}{2}}
\sum_{\pi\in \mathcal{S}_n} (-s)^{-l(\pi)}w_\pi\ .$$
 Here $l(\pi)$ is the length of the permutation $\pi$.
 We work in the generic case.
This means that, in the domain $k$,
the quantum integers $[j]=\frac{s^j-s^{-j}}{s-s^{-1}}$ are 
asked to be invertible
for every $j>0$.


 For a Young diagram $\lambda$ of size $n$, we denote
by $\square_\lambda$ the object of the category ${{H}}$ formed
with one point  for each cell $c$ of $\lambda$, 
equipped with the standard trivialization. If $c$ has index $(i,j)$
($i$-th row, and $j$-th column),
then the corresponding point in ${D}^2$ is
$\frac{j+i\sqrt{-1}}{n+1}$.
Following Aiston and Morton \cite{AM}, we can define in
${H}_{\Box_\l}$ a minimal idempotent $y_\l$
which is a version of the corresponding Young idempotent
of the symmetric group algebra.
 The idea of the construction is to insert symmetrizers
along rows and antisymmetrizers along columns, and then to normalize.
 A skein computation of the normalizing coefficient
(of slightly different idempotents) appeared in \cite{Yo}.
 Details about the construction of these idempotents can also be found
in \cite{Bhec}.

 A standard tableau $t$ with shape a Young diagram
 $\lambda=\lambda(t)$ is a labeling of the cells, with the integers
 $1$ to $n$, which is increasing along rows and columns.
 We denote by $t'$ the tableau obtained by removing
 the cell numbered by $n$.
 We define $\alpha_t\in {{H}}(n,\square_\lambda)$ and
$\beta_t\in {{H}}(\square_\lambda,n)$
 by
 $$\alpha_1=\beta_1=\oneun\ ,$$
 $$\alpha_t=(\alpha_{t'}\otimes \oneun)\varrho_t y_\lambda\ ,$$
 $$\beta_t=y_\lambda \varrho_t^{-1}(\beta_{t'}\otimes\oneun)\ . $$
 Here $\varrho_t\in {{H}}(\square_{\lambda(t')}\otimes 1,
\square_\lambda)$ is a standard isomorphism.

The following theorems are shown in \cite{Bhec}.
\begin{thm}\label{astruct}
The family $\alpha_t\beta_\tau$ for all standard tableaus $t,\tau$
such that
$\lambda(t)=\lambda(\tau)$
 forms a basis for ${{H}}_n$,
and (here $\delta_{\tau s}$ is the Kronecker delta)
$$\alpha_t\beta_\tau\alpha_s\beta_\sigma=
\delta_{\tau s} {\alpha_t\beta_\sigma}\ .$$
\end{thm}
This gives explicitly  an algebra isomorphism
$$\bigoplus_{|\lambda|=n} \mathcal{M}_{d_\lambda}(k)\approx {H}_n
\ ,$$
where $d_\lambda$ is the number of standard tableaus with shape $\lambda$,
and $\mathcal{M}_{d_\lambda}(k)$ is the algebra of $d_\lambda\times d_\lambda$
matrices with coefficients in $k$.
The diagonal elements $p_t=\alpha_t\beta_t$
are the {\em path idempotents} described in
 \cite{Wenzl}.
The {\em minimal central idempotent} corresponding to the partition
$\l$ is $$z_\l=\sum_{\l(t)=\l}p_t\ .$$
The minimal idempotents $y_\l$ and the path idempotents $p_t$
satisfy the following branching formula.
\begin{thm}[Branching formula]\label{abranch}
$$\begin{array}{rcl}y_\lambda\otimes \oneun&=&
\displaystyle \sum_{\genfrac{}{}{0pt}{2}{\lambda\subset \mu}{|\mu|=|\lambda|+1}}
(y_\lambda\otimes \oneun)y_\mu (y_\lambda\otimes \oneun)\ ,\\
\rule{0pt}{18pt}p_t\otimes \oneun&=&
\displaystyle\sum_{\tau'=t}
p_\tau\ .\end{array}$$
\end{thm}
\noindent We have omitted in these  formulas the standard
isomorphisms respectively between $\square_\lambda\otimes 1$ and $\square_\mu$,
and between $(n-1)\otimes 1$ and $n$.

The result for the quantum dimensions
is given in the following theorem \cite{Wenzl}.
Here we denote by $\la\ \ra^h$ the quantum trace in the
ribbon Hecke category.

\begin{thm}[Quantum dimension]\label{aqdim}
$$ \langle y_\lambda\rangle^h=
\prod_{\mathrm{cells}}{\frac{\alpha s^{cn(c)}-\alpha^{-1}s^{-cn(c)}}
{s^{hl(c)}-s^{-hl(c)}}}$$
\end{thm}

The assertion above can  be proven by  a skein calculation
(see  \cite[Prop.2.4]{Yo} or \cite{Aiston2}). 
 Here is a sketch of the proof.
We first check the formula for columns $1^n$, by using
the recursive formula for the antisymmetrizers $y_{1^n}$.
We note \cite[Ch.1]{Macdo} that the right hand side in Theorem \ref{aqdim}
is the Schur polynomial in the $\la y_{1^n}\ra^h $.
We then proceed recursively on the number of cells.

If $\l$ contains two distinct sub-diagrams
$\m$ and $\nu$ with
$|\mu|=|\nu|=|\lambda|-1$, then we get the result by considering
$(y_\m\otimes \one\!{}_1)(y_\n\otimes \one\!{}_1)$
(considered as a composition of morphisms
in $H_{\Box_\l}$).
 This defines a quasi-idempotent which can be normalized.
Whence  we get  a 
minimal idempotent which belongs to the simple component
indexed by $\l$, whose quantum trace gives the required formula.

 We obtain the remaining cases, namely the rectangular diagrams,
$\l=(j,\dots,j)$ by using the branching
formula for $\mu=(j,\dots,j,j-1)$.


\section{Idempotents in Birman-Murakami-Wenzl algebras}\label{s_idem}
The quotient of $K_n$ by the ideal $I_n$ generated by $h_{n-1}$
is isomorphic to the Hecke algebra ${H_n}$.
 We denote by $\pi_n$
the canonical projection map
$$\pi_n:K_n 
{\longrightarrow} H_n\ .$$
The main idea of our construction is to define
a multiplicative section $s_n: H_n \to K_n$ and to use it 
for the transport of the
idempotents from the Hecke algebra to the Birman-Murakami-Wenzl one.
 In this section, we suppose that the ground ring $k$
is the field $\bb{Q}(\a,s)$. The computation of the quantum dimensions
in Section \ref{sectionqdim} will permit to discuss the non generic 
case.

\begin{thm}\label{tr}
There exists a unique multiplicative homomorphism
$s_n: H_n\to K_n$, such that 
 $$\pi_n\circ s_n= \mathrm{id}_{H_n}\text{ and }$$
 $$s_n(x)y=ys_n(x)=0 \;\;{\it for}\;\; \forall x\in H_n\;{\it and}\;
 \forall y\in I_n$$
\end{thm}
\begin{cor}\label{cor1}
$K_n\cong I_n\oplus H_n\cong I_n\oplus  
\left(\oplus_{|\l|=n} {\cal M}_{d_\l} (k)\right) $. 
\end{cor}
\noindent
The theorem above gives minimal central
idempotents in $K_n$,
$$\tilde z_\l=s_n(z_\l)\ ,$$ 
and also minimal path  idempotents, $$\tilde p_t=s_n(p_t)\ .$$
 The quantum trace of $\tilde p_t$ depends only on $\l=\l(t)$;
it is denoted by $\la\l\ra$.
 As we did in Section \ref{hecke}, for each Young diagram
$\l$, we consider an object $\Box_\l$ whose points
correspond to the cells of $\l$.
From Theorem \ref{tr}, we get the section 
$s_{\Box_\l}:H_{\Box_\l}\rightarrow K_{\Box_\l}$,
and we obtain a  minimal idempotent in $K_{\Box_\l}$,
$$\tilde y_\l=s_{\Box_\l}(y_\l)\in K_{\Box_\l}\ .$$

\begin{lem}\label{orth}
 If $\l$ and $\mu$ are two distinct Young diagrams with
the same size $|\l|=|\mu|=n$, then for every
$x\in K(\Box_\l,\Box_\m)$ one has
$\tilde y_\l x\tilde y_\m=0$ .\\

\end{lem}
\begin{proof}
We have  $x=\tilde h+y$ with $\tilde h=s_n(h)$, $h\in H_n$ and
$y\in I_n$. We omit here the isomorphisms between 
$K(n,\Box_\l)$ and $K(n,\Box_\m)$. Then 
$$s_n(y_\l)x s_n(y_\m)=s_n(y_\l)\tilde h s_n(y_\m)=s_n(y_\l h y_\m) .$$
The result follows from the corresponding 
property in the category $H$.
\end{proof}
 We will need the following absorbing property
which also results from computation in the Hecke algebra
\cite[Cor.1.10]{Bhec}.
\begin{lem} If $\nu$ is a Young diagram obtained from the Young diagram 
$\lambda$ by removing one cell, then one has:
$$\label{abs}\tilde y_\l(\tilde y_\nu\otimes\oneun)\tilde y_\l=\tilde y_\l\ .$$
\end{lem}
\v8
\noindent
{\it Proof of Theorem \ref{tr}.}
If the section $s_n$ exists, then it is unique. This can be seen as follows.
Let $U_n$ be the central idempotent corresponding to the factor $I_n$. 
If $s'_n$ is another section, then we have for every $x\in H_n$
$$s_n(x)-s'_n(x)=(s_n(x)-s'_n(x))\, U_n=0\ .$$

We will construct the section $s_n$ by induction on $n$.
The result is certainly true for $n=1$, since
we have that $K_1\approx H_1\approx k$.


Let us assume that we have constructed $s_{m}$
satisfying the conditions of the theo\-rem, for every $m<n$,
 so that we have
 minimal idempotents, $\tilde p_t$,
for every $t$ such that $|\l(t) |<n$.

 Let $\l$ be a Young diagram whose size is $|\l|=n-1$, then we have
a minimal idempotent $\tilde y_\l\in K_{\Box_\l}$.
 If $\nu$ is a Young diagram included in $\l$, such that
$|\n|=n-2$, then we define 
$\tilde y_{(\l,\nu)}\in K_{\Box_\l\otimes 1}$, by
\be \label{yln}\tilde y_{(\l,\nu)}=\frac{\la \n\ra}{\la\l\ra}
(\tilde y_\l\otimes \one\!{}_1)(\tilde y_\nu\otimes h_{1})
 (\tilde y_\lambda\otimes \one\!{}_1)\ .\ee
Here, the standard isomorphisms between
$\Box_\l\otimes 1$ and $\Box_\nu\otimes 2$ are omitted.
We need that the quantum dimensions $\la\l\ra$, with $|\l|=n-1$
are not zero. This result  follows from \cite[Theorem 3.7]{BW}, and will
 also be proved, by considering the specializations
corresponding to Brauer algebras in the next section.

\begin{lem}\label{ortho} a) If $\nu$ and $\mu$ are two distinct
Young diagrams 
of  size $n-2$, included in $\l$, then
$\tilde y_{(\l,\nu)}\tilde y_{(\l,\mu)}=0\ .$\\
b) If $\nu$ is a
Young diagram 
of  size $n-2$, included in $\l$, then $\tilde y_{(\l,\nu)}$
is an idempotent.
\end{lem}
\begin{proof}
By the induction hypothesis we can  apply Lemmas \ref{orth}
and \ref{abs} to Young diagrams of size $m<n$.
The statement $a)$ 
 follows then from Lemma \ref{orth}
applied to $\n$ and $\m$ with $|\n|=|\m|=n-2$.

The square of $\tilde y_\nu$ is equal to
$\left(\frac{\langle\nu\rangle}{\langle\l\rangle}\right)^2$ times
the skein element represented by the following tangle.
$$\psdiag{15}{45}{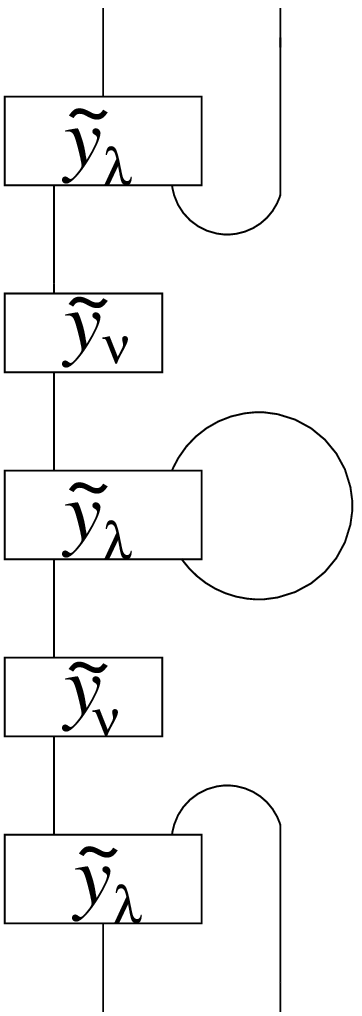}\; $$
 Let us consider the intermediate morphism
$$(\tilde y_\nu \otimes\cup_{\l/\nu})
(\tilde y_\l \otimes\oneun)(\tilde y_\nu \otimes\cap_{\l/\nu})$$
(the subscript in $\cup_{\l/\nu}$ and $\cap_{\l/\nu}$
indicate which isomorphism in $K(\Box_\l\otimes 1,\Box_\n\otimes 2)$
is used).
The  minimality  of the idempotent
$\tilde y_\nu$ implies that this  morphism is equal to $\tilde y_\nu$,
up to a coefficient which is obtained by considering the trace (we use
the absorbing property \ref{abs}).
\be\label{ptrace}(\tilde y_\nu \otimes\cup_{\l/\nu})
(\tilde y_\l \otimes\oneun)(\tilde y_\nu \otimes\cap_{\l/\nu})
=\frac{\la\l\ra}{\la\nu\ra}\tilde y_\nu\ .\ee
  Statement $b)$ follows.
\end{proof}
 
Let $t$ be a standard tableau whose size is $n-1$ , with shape
 $\lambda(t)=\lambda$. 
 We define $\mathfrak{a}_t\in {K}(n-1,\square_\lambda)$ and
$\mathfrak{b}_t\in {K}(\square_\lambda,n-1)$
by lifting to the category $K$ the elements
$\alpha_t$ and $\beta_t$ defined in Section \ref{hecke}.
If $\alpha_t$ and $\beta_t$ are represented by linear combinations of `braids'
(elements of the braid groupoid),
then $\mathfrak{a}_t$ and $\mathfrak{b}_t$ are given by 
the following formulas.
$$\mathfrak{a}_t=s_{n-1}(\one\!_{n-1})\a_t=\a_t\tilde y_\l\ ,$$
$$\mathfrak{b}_t=\beta_t s_{n-1}(\one\!_{n-1})=\tilde y_\l \beta_t\ .$$
We then have the formula
$$\tilde p_t=s_{n-1}(p_t)=\mathfrak{a}_t\mathfrak{b}_t\ .$$
For a  Young diagram $\nu$ with
$n-2$ cells,
included in $\l=\l(t)$, we define 
$\tilde p_{(t,\nu)}$ and $\tilde{p}^+_t\in K_n$, by
$$\tilde p_{(t,\nu)}=(\mathfrak{a}_t\otimes \one\!_1)\tilde y_{(\l,\nu)}
(\mathfrak{b}_t\otimes \one\!_1)\ ,$$
$$\tilde{p}^+_t= \tilde{p}_t\otimes \one\!{}_1 - 
\sum_{\genfrac{}{}{0pt}{2}{\n\subset \l(t)}
{|\n|=n-2}} \ \tilde p_{(t,\nu)}
\ .$$
Using  Lemma \ref{ortho}, we get the following.
\begin{lem}\label{ortho2}
i) $\tilde p_{(t,\nu)}\tilde p_{(\tau,\mu)}=\delta_{t\tau}\delta_{\nu\mu}\tilde p_{(t,\nu)}$; $\;\;\;$
ii) $ \tilde{p}^+_t \tilde{p}^+_\tau =
 \d_{t \tau} \tilde{p}^+_t$ . 
\end{lem}
We  define a linear homomorphism $s'_n$
from the braid group algebra $k[B_n]$ to $K_n$ by 
\be \label{defsn}\forall x\in k[B_n]\;\;\;\;
 s'_n(x)= \sum_{|\l(t)|=|\l(\tau)|=n-1} 
\tilde p^+_t  \;x \; \tilde p^+_\tau\ .\ee

We will show that
this homomorphism induces a well defined linear map
$s_n: H_n\rightarrow K_n$ which is a section of $\pi_n$, and  
prove multiplicativity.
The proof will be complete with the two following lemmas.
$\hfill\Box$

\v8
\begin{lem}\label{proofa}
$ \tilde{p}^+_t  y =y \tilde{p}^+_t =  0 \;\;{\it for}\;\; \forall y\in I_n$. 
\end{lem}

\begin{proof}
 We want to show
that
$\tilde{p}^+_t  y = 0$  for any $y\in I_n$.
 We write  $$y=(a\otimes \oneun) h_{n-1}(b\otimes \oneun) \ .$$
 By the induction hypothesis, we have the result if $a$
is in $I_{n-1}$. So it is enough to consider the case where
$a=s_{n-1}(x)$ for $x$ an element
of the matrix units basis  described
in Section \ref{hecke}, i.e. for $a=\mathfrak{a}_\sigma
 \mathfrak{b}_\tau$,
where $\tau$ and $\sigma$ are standard tableaus  with the same shape,
whose size
is $n-1$.
 If $\sigma\neq t$, then we have
 $\tilde{p}_t\mathfrak{a}_\sigma=0$ and the result follows.
It remains to check the case where $\sigma=t$.
$$\begin{array}{rcl}
\tilde{p}^+_t (\mathfrak{a}_t
 \mathfrak{b}_\tau\otimes \oneun)h_{n-1}&=&
\tilde{p}^+_t (\mathfrak{a}_t
 \otimes \oneun)(\mathfrak{b}_{\tau'}\otimes h_{1})
\\&=&
(\tilde{p}_t\otimes \one\!{}_1) (\mathfrak{a}_t
 \otimes \oneun)(\mathfrak{b}_{\tau'}\otimes h_{1})
\\&&
{\displaystyle-\sum_{\genfrac{}{}{0pt}{2}{\nu\subset \l}
{|\n|=n-2}} \frac{\la\nu\ra}{\la\l(t)\ra}
(\mathfrak{a}_t\otimes\oneun)(\tilde y_\nu\otimes h_{1})
(\mathfrak{b}_t\otimes\oneun)(\mathfrak{a}_t
 \otimes \oneun)(\mathfrak{b}_{\tau'}\otimes h_{1})}
\\&=&
(\mathfrak{a}_t
 \otimes \oneun)(\mathfrak{b}_{\tau'}\otimes h_{1})
\\&&
{\displaystyle- \frac{\la\l(\tau')\ra}{\la\l(t)\ra}
(\mathfrak{a}_t\otimes\oneun)(\tilde y_{\l(\tau')}\otimes h_{1})
(\tilde y_{\l(t)}
 \otimes \oneun)(\mathfrak{b}_{\tau'}\otimes h_{1})}\\
&=&0
\end{array}   
$$ 
The result $\tilde{p}^+_t  y = 0$ follows; 
$y\tilde{p}^+_t  = 0$ can be  obtained similarly.


\end{proof}
\begin{lem}
The map $s'_n$ induces a well defined multiplicative homomorphism
$$s_n: H_n\rightarrow K_n$$
such that $\pi_n\circ s_n=\mathrm{id}_{H_n}$.
\end{lem}
\begin{proof}

 From Lemma \ref{proofa}, we can see that the Homfly skein relation
is respected, whence we have that $s_n$ is well defined.
We have that $\pi_n(\tilde p_t^+\otimes \oneun)=p_t\otimes \oneun$.
 This implies that $\pi_n\circ s_n=\mathrm{id}_{H_n}$.
The computation below shows the multiplicativity.
$$\begin{array}{rcl}
 s_n(x) s_n(y)&=&
\displaystyle 
\sum_{t,\sigma,\tau} \tilde p^+_t x \tilde p^+_\sigma  y \tilde p^+_\tau
\\&=&\displaystyle 
\sum_{t,\sigma,\tau} \tilde p^+_t 
x (\tilde p_\sigma\otimes \oneun)  y 
\tilde p^+_\tau
\\&=&\displaystyle 
 s_n\left(\sum_\sigma
x ( p_\sigma\otimes \oneun)  y \right)
\\&=& \displaystyle 
 s_n(xy)\ .
\end{array}$$
\end{proof}

\section{Brauer algebras}\label{brauer}
Brauer centralizer  algebras were introduced in \cite{Br}
in relation with the representation theory of the orthogonal and
symplectic groups (see also \cite{Weyl}).
Their structure was obtained by Wenzl in \cite{WX}.
 We emphasize also Nazarov's computations of the action of  generators
on their irreducible
representations in \cite{Naz}. His work includes the dimension
formulas and inspired our computation of the quantum dimensions in
Section \ref{sectionqdim}.

We have defined the Birman-Murakami-Wenzl algebras $K_n$ by using
Kauffman skein theory. Brauer algebras can be defined
in a similar way by using the classical version of 
Kauffman relations given in Figure \ref{clkaufrel}.
\begin{figure}[h]
$$
\psdiag{3}{9}{px}\;=\;\psdiag{3}{9}{nx} 
$$
$$
\psdiag{3}{9}{ptwist}\;=\;\; \psdiag{3}{9}{1}\;$$
$$L\;\amalg \;\bigcirc\; = \;N
\;\, L  $$
\caption{\label{clkaufrel} Classical Kauffman skein relations}
\end{figure}

Here any coefficient ring is allowed,
and $N$ could be an indeterminate.
 If $N$ is a natural number, then we obtain Brauer algebras
with complex coefficients 
as a specialization of Birman-Murakami-Wenzl algebras,
with coefficient ring $\bb{C}[s^{\pm1}]$ and 
 $\a=s^{N-1}$, by setting $s=1$.
 It is a classical fact \cite[Ch.5]{Weyl}, that there exists an algebra
homomorphism $\Phi_n$ from this Brauer algebra, denoted
by $D_n(N)$, to the centralizer algebra
$End_{O(N)}(V^{\otimes n})$, where $V=\bb{C}^N$ is the 
fundamental representation of the orthogonal group $O(N)$.
 If  we denote by $(u_1,\dots, u_N)$ the canonical basis of $V$,
then $\Phi_n$ is defined on the generators by
$$\Phi_n(e_i).u_{j_1}\otimes\dots u_{j_i}\otimes u_{j_{i+1}}\otimes\dots u_{j_n}=
u_{j_1}\otimes\dots  u_{j_{i+1}}\otimes u_{j_i}\otimes\dots u_{j_n}\ ,$$
$$\Phi_n(h_i).u_{j_1}\otimes\dots u_{j_i}\otimes u_{j_{i+1}}\otimes\dots u_{j_n}=
\delta_{{j_i}{j_{i+1}}}\sum_{\nu=1}^n u_{j_1}\otimes\dots  u_{j_\nu}\otimes 
u_{j_\nu}\otimes\dots u_{j_n}\ .$$
In the above, $\delta$ is the Kronecker delta.
 For $N\geq n$, this homomorphism is injective.
 In fact the above assignment extends to a monoidal functor
from the specialized BMW category, to the linear category.
  This functor is compatible with the (trivial)
ribbon structures on these categories, and so it respects the `quantum' traces.
 In the case of the linear category the quantum trace coincides with the usual 
one. 
 We get that,
with the given specializations,
$$\forall x\in D_n(N)\;\;\;\; \langle x\rangle = 
\mathrm{trace}(\Phi_n(x))\ .$$
In particular, if $x$ is a non trivial idempotent,
then its quantum dimension is a non zero natural number. 

The quotient of $D_n(N)$ by the ideal $I_n$ generated by $h_{n-1}$
is isomorphic to the symmetric group algebra $\bb{C}[\cal{S}_n]$,
the classical counterpart of the Hecke algebra.
 We denote by $\pi_n$
the canonical projection map
$$\pi_n:D_n(N)\;
{\longrightarrow}\; \bb{C}[\cal{S}_n]\ .$$

\begin{thm}\label{cltr}
If $N$ is  an integer greater or equal to $n$,
then, there exists a multiplicative homomorphism
$s_n: \bb{C}[\cal{S}_n]\to D_n(N)$, such that 
 $$\pi_n\circ s_n= \mathrm{id}_{\bb{C}[\cal{S}_n]}\text{ and }$$
  $$\forall x\in \bb{C}[\cal{S}_n]\;\;\; 
\forall y\in I_n\;\;\; s_n(x)y=ys_n(x)=0 \ .$$
\end{thm}
\begin{rem}\label{r1}
As a corollary we get a non trivial minimal idempotent
$\tilde p_t=s_n(p_t)$ for every standard tableau of size $n$,
whose quantum trace is a non zero natural number.
 Here $p_t$ is the  minimal path idempotent in the symmetric group
algebra.
\end{rem}
\begin{proof}
The recursive construction of the preceding section can be done.
The only point to check is that at each step the quantum dimensions
$\langle\l\rangle$ are not zero.
 By the induction hypothesis, we have non trivial minimal
idempotents $\tilde p_t=s_n(p_t)$ for every standard tableau of size $n-1$.
If the shape of $t$ is $\l$, then we have that
$\langle\l\rangle=\langle \tilde p_t\rangle=\mathrm{trace}
(\Phi(\tilde p_t))$
is a non zero integer.
\end{proof}
\begin{rem}
The above completes the proof of Theorem \ref{tr}. At each step
of the recursive construction, we needed that the quantum
dimensions $\langle\l\rangle$ with $|\l|=n-1$ are non zero.
This is the case, because they become non zero integers
if we apply the rank $N\geq n$ Brauer specialization.
 The same remark shows that Theorem \ref{cltr} holds if $N$ is generic
(e.g. for the Brauer algebra with ground field $\bb{Q}(N)$).
\end{rem}

\section{Matrix units in Birman-Murakami-Wenzl algebras}\label{semi}
In this section we will describe a matrix units basis in $K_n$
(compare with \cite{RW}),
and show the branching formula.
 Recall that in the categories $H$ and $K$,
for each Young diagram
$\l$, we have defined an object $\Box_\l$ whose points
correspond to the cells of $\l$.
 In the proof of Theorem \ref{tr}, we have used the
section  $s_{\Box_\l}:H_{\Box_\l}\rightarrow K_{\Box_\l}$ to define
 a minimal idempotent
$\tilde y_\l=s_{\Box_\l}(y_\l)\in K_{\Box_\l}$.

A sequence
$\Lambda=(\Lambda_1,\dots,\Lambda_{n})$ of Young diagrams, in which
two consecutive diagrams $\Lambda_i$ and $\Lambda_{i+1}$
differ by exactly one cell will be called an  up and down tableau
of length $n$, and shape $\Lambda_n$.
 Those up and down tableaus of length
$n$  such that $n=|\Lambda_{n}|$ (up tableaus) correspond bijectively with
standard tableaus as described in Section \ref{hecke}.

 For an up and down tableau $\Lambda$
of length $n$, we denote by $\Lambda'$ the tableau of length
$n-1$ obtained by removing
 the last Young diagram in the sequence $\Lambda$.
 We define $\mathfrak{a}_\Lambda\in {{K}}(n,\square_\Lambda)$ and
$\mathfrak{b}_\Lambda\in {{K}}(\square_\Lambda,n)$
 by
 $$\mathfrak{a}_1=\mathfrak{b}_1=\oneun\ ,$$
if $|\Lambda_n|=|\Lambda_{n-1}|+1$, then
 $$\mathfrak{a}_\Lambda=(\mathfrak{a}_{\Lambda'}\otimes \oneun) \tilde y_{\Lambda_n}\ $$
 $$\mathfrak{b}_\Lambda=\tilde y_{\Lambda_n} (\mathfrak{b}_{\Lambda'}\otimes\oneun)\ , $$
if $|\Lambda_n|=|\Lambda_{n-1}|-1$, then
$$\mathfrak{a}_\Lambda=\frac{\la\Lambda_n\ra}{\la\Lambda_{n-1}\ra}
(\mathfrak{a}_{\Lambda'}\otimes \oneun) 
(\tilde y_{\Lambda_n}\otimes\cap)\ $$
 $$\mathfrak{b}_\Lambda=(\tilde y_{\Lambda_n}\otimes\cup)
 (\mathfrak{b}_{\Lambda'}\otimes\oneun)\ . $$

Here we have omitted the standard isomorphism
 in ${{K}}(\square_{\Lambda_{n-1}}\otimes 1,
\square_{\Lambda_n})$ and in ${{K}}(\square_{\Lambda_{n-1}}\otimes 1,
\square_{\Lambda_n}\otimes 2)$. Note that for an up tableau, the definition is
coherent with the one given for the corresponding standard tableau in the proof
of Theorem \ref{tr}.

\begin{thm}\label{kbase}
a) The family $\mathfrak{a}_\Lambda\mathfrak{b}_\Xi$ for all up and down 
 tableaus $\Lambda$,$\Xi$ of length $n$,
such that
$\Lambda_n=\Xi_n$ forms a basis for $K_n$,
and 
$$\mathfrak{a}_\Lambda\mathfrak{b}_\Xi\mathfrak{a}_L\mathfrak{b}_X=
\delta_{\Xi L}\ {\mathfrak{a}_\Lambda\mathfrak{b}_X}\ .$$
b) There exists an algebra isomorphism
$$\bigoplus_{|\lambda|=n,n-2,\dots} 
\mathcal{M}_{d_\lambda^{(n)}}(k)\approx {K}_n
\ ,$$
where $d_\lambda^{(n)}$ is the number of the  up and down
 tableaus of length $n$ with shape $\lambda$,
and $\mathcal{M}_{d}(k)$ is the algebra of $d\times d$
matrices with coefficients in $k$.
\end{thm}

The diagonal elements $q_\Lambda=\mathfrak{a}_\Lambda\mathfrak{b}_\Lambda$
are the path idempotents associated with the inclusions $K_i\subset K_{i+1}$
the minimal central idempotent corresponding to the partition
$\l$  is 
$$z^{(n)}_\l=\sum_{\Lambda_n=\l}q_\Lambda\ .$$
If $\Lambda$ corresponds to a standard tableau $t$ (Hecke part),
then one has $q_\Lambda=\tilde p_t$, and if 
$|\lambda|=n$ then $z^{(n)}_\l=\tilde z_\lambda$.
\begin{proof}
From Lemma \ref{orth} if $|\l|=|\m|$,
and from Lemma \ref{ortho2} if $|\l|\neq|\m|$, we get
\be\tilde y_\l K(\Box_\l,\Box_\m)\tilde y_\m=0 \ .\ee
This implies that $\mathfrak{b}_\Lambda\mathfrak{a}_\Xi=0$ if $\Lambda\neq\Xi$.

 Using the  properties \ref{abs} and \ref{ptrace}, we show that, 
if $n$ is the length of $\Lambda$, we have
$$\mathfrak{b}_\Lambda\mathfrak{a}_\Lambda=\tilde y_{\Lambda_n}\ .$$
Hence we have that
$$\mathfrak{a}_\Lambda\mathfrak{b}_\Xi\mathfrak{a}_L\mathfrak{b}_X=
\delta_{\Xi L}\ {\mathfrak{a}_\Lambda\mathfrak{b}_X}\ .$$
The independence follows.

To show that the family $\mathfrak{a}_\Lambda\mathfrak{b}_\Xi$ generate $K_n$,
we proceed recursively on $n$.
Using the map $s_n$, we see from the known result in $H_n$
that the $\mathfrak{a}_\Lambda\mathfrak{b}_\Xi$, where $\Lambda_n$ and $\Xi_n$ 
is the same  Young diagram with $n$ cells,
generate the Hecke part $\widetilde H_n=s_n(H_n)$
of $K_n$. It remains to consider $I_n$. From the induction
hypothesis, we get that $I_n$ is generated
by the following elements
$$(\mathfrak{a}_\Lambda\otimes \oneun)(\mathfrak{b}_\Xi\otimes \oneun)
h_{n-1}(\mathfrak{a}_L\otimes \oneun)(\mathfrak{b}_X\otimes \oneun)\ .
$$
These are zero if $\Xi'\neq L'$, and  else
are equal to
$$\frac{\langle \Lambda_{n-1}\rangle}
{\langle\l\rangle}\mathfrak{a}_{(\Lambda,\l)}\mathfrak{b}_{(X,\l)}\ ,$$
where $\l=\Xi'_{n-2}=L'_{n-2}$. 
\end{proof}
The path idempotents $q_\Lambda$ satisfy the following branching formula.
\begin{thm}[Branching formula]\label{branch}
$$ q_\Lambda\otimes\oneun=\sum_{\Xi'=\Lambda}q_\Xi\ .$$
\end{thm}
We recall that in the above formula
the shape of $\Xi$ is either one cell more, either one cell less
than the shape of $\Lambda$. 

\begin{proof}
The coordinates of $ q_\Lambda\otimes\oneun$ in the standard basis
are obtained by computing
$\mathfrak{b}_\Xi(q_\Lambda\otimes\oneun)\mathfrak{a}_L$.
The result is zero unless $\Xi$
and $L$ have the same shape $\mu$
and $L'=\Xi'=\Lambda$; and in the latter case
the result is $\tilde y_\mu$.
\end{proof}

As a corollary we also have a branching formula
for the minimal idempotents $\tilde y_\l$ (the obvious isomorphisms
are omitted).
\begin{cor}\label{branch2}
$$\tilde y_\l\otimes \one\!{}_1=
\sum_{\genfrac{}{}{0pt}{2}{\l\subset \mu}{|\mu|=|\l|+1}}
(\tilde y_\lambda\otimes \oneun)\tilde y_\mu (\tilde y_\lambda\otimes \oneun)
+
\sum_{\genfrac{}{}{0pt}{2}{\nu\subset \lambda}{|\nu|=n-1}}
\tilde y_{(\lambda,\nu)}\ .
$$
\end{cor}
\begin{proof}
We can decompose
$\tilde y_\l\otimes 1$ by using the minimal central idempotents
in $K_{|\l|+1}$. We get
$$
\tilde y_\l\otimes 1=
\sum_\Xi (\tilde y_\l\otimes \oneun)q_\Xi 
(\tilde y_\l\otimes \oneun)\ .$$
In the above, only those up and down tableaus $\Xi$ with
$\Xi_{|\l|}=\l$ contribute. 
\end{proof}

\section{Braiding and twist coefficients}

\begin{pro}
i) (Twist coefficient)
Let $\tilde{y}_\m\in K_{\Box_\m}$  be the minimal idempotent, then
\be\label{fr}\psdiag{6}{18}{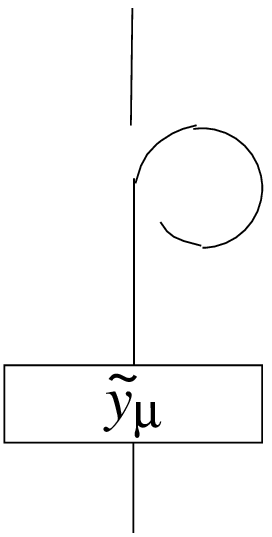}= \a^{|\m|} s^{2\sum_{c\subset \m} cn(c)} \tilde y_\m\; .\ee
Here the content of a cell $c$ in the $i$-th row and 
the $j$-th column of $\mu$
is defined by $cn(c):=j-i$.

ii) \label{brad}(Braiding coefficient)
Suppose that $\l \subset \m$ and $\m-\l$ contains only one cell $c$.
 Then
\be\label{oi}
\psdiag{15}{45}{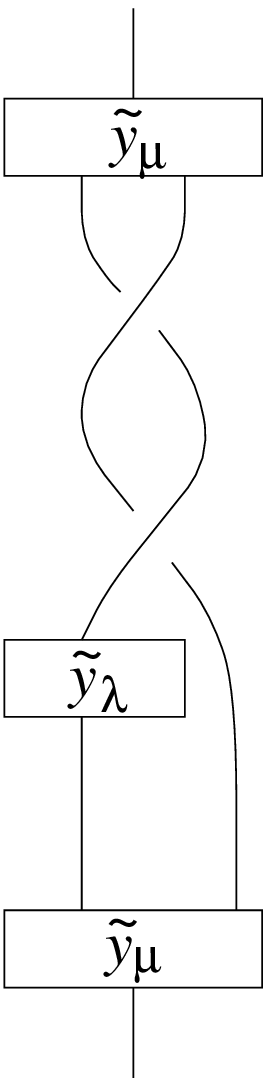}= s^{2cn(c)} \tilde{y}_{\m}\, .\ee
Suppose that $\m \subset \l$ and $\l-\m$ contains only one cell $c$, then
\be\label{io}
\psdiag{15}{45}{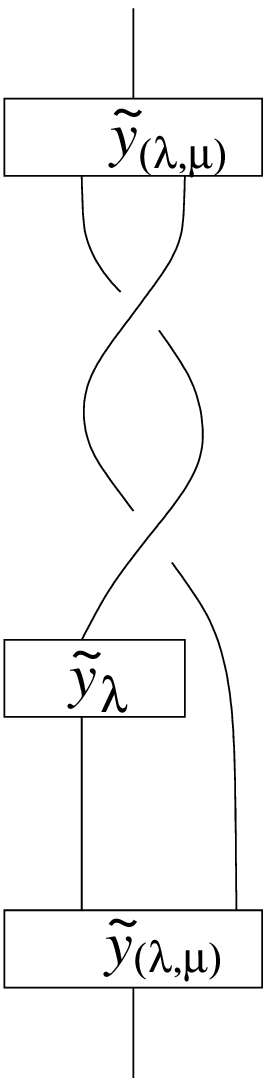}= \a^{-2} s^{-2cn(c)} \tilde{y}_{(\l,\m)}\, .\ee
\end{pro}
\begin{proof}
The statements (\ref{fr}) and (\ref{oi}) follow  from 
the corresponding in the Hecke algebra
(see  \cite[Prop.1.11]{Bhec}).
Using the  definition of the idempotent $\tilde y_{(\l,\mu)}$, we can bring (\ref{io})
to the form
$$\left(\frac{\la\m\ra}{\la \l\ra}\right)^2
(\tilde y_\l\otimes\one\!{}_1)\left((\tilde y_\m x \tilde y_\m)
\otimes h_{\l/\m}\right)(\tilde y_\l\otimes\one\!{}_1)\ ,$$
where $x\in K_{\Box_\m}$ is depicted below.
$$\psdiag{12}{36}{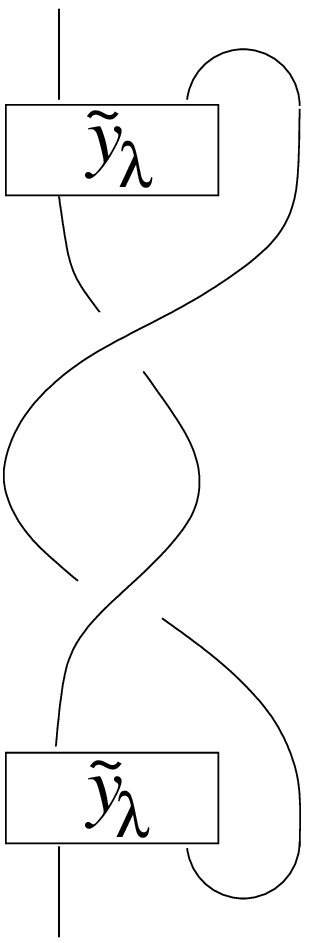}$$
By the Schur lemma, $\tilde y_\m x \tilde y_\m= c \tilde y_{\m}$ with $c\in k$.
Taking the quantum trace of this morphism, we get
$$\la \tilde y_\m x \tilde y_\m\ra
= \a^{-2}s^{-2 cn(c)} \la \l\ra = c \la \m \ra\; .$$ 
The first equality is due to the Kauffman skein relations and 
the statement (\ref{oi}) above.
\end{proof}

\section{Quantum dimension}\label{sectionqdim}
The formula for the quantum dimension $\la \l\ra$ was obtained by Wenzl
\cite[Theorem 5.5]{Wbcd}. The proof there rests on the representation theory
of the quantum group $U_qso(2n+1)$. We give below an 
alternative  proof 
for this formula. Our method is inspired by the Nazarov computation of
the matrix elements of  the action of the hook generators on
the canonical basis of the irreducible representations of Brauer algebras 
\cite{Naz}. Here we work with scalar field $k=\Q(s,\a)$, hence we have that
all the quantum dimensions $\langle \l\rangle$ are invertible.

Suppose that $\l=(\l_1,\dots,\l_m)$ 
is obtained from $\m$ by adding one cell in the
$i$th row. Let $l$ be the number of pairwise distinct
rows in the diagram $\mu$. Then one can obtain $l+1$ diagrams by 
adding a cell to $\mu$, and $l$ diagrams by 
removing a cell from $\mu$. Let $c_1,\dots,c_{l+1}$ and
$d_1,\dots,d_{l}$ be the contents  (defined in Prop.\ref{brad}) 
of these cells respectively.
 Denote by $b_1,\dots,b_{2l+1}$ the scalars
$$\a s^{2c_1},\dots,\a s^{2c_{l+1}},$$
$$\a^{-1}s^{-2d_1}\dots,\a^{-1}s^{-2d_{l}} ,$$
and by $b$ the value among $b_1,\dots,b_{l+1}$
corresponding to the diagram $\l$.
\begin{thm}\label{recdim}
One has
\be\label{rec}\frac{\la \l\ra}{\la \m\ra}= 
\alpha b^{-1}\left(\frac{b-b^{-1}}{s-s^{-1}}+1\right)
\prod_{b_j\neq b}\frac{b-b_j^{-1}}{b-b_j}\ .\ee
\end{thm}
\begin{proof}
Let $\tau_n$ be the element of the algebra $K_n$
defined below.

$$\label{taun}
\tau_n=\ \psdiag{15}{45}{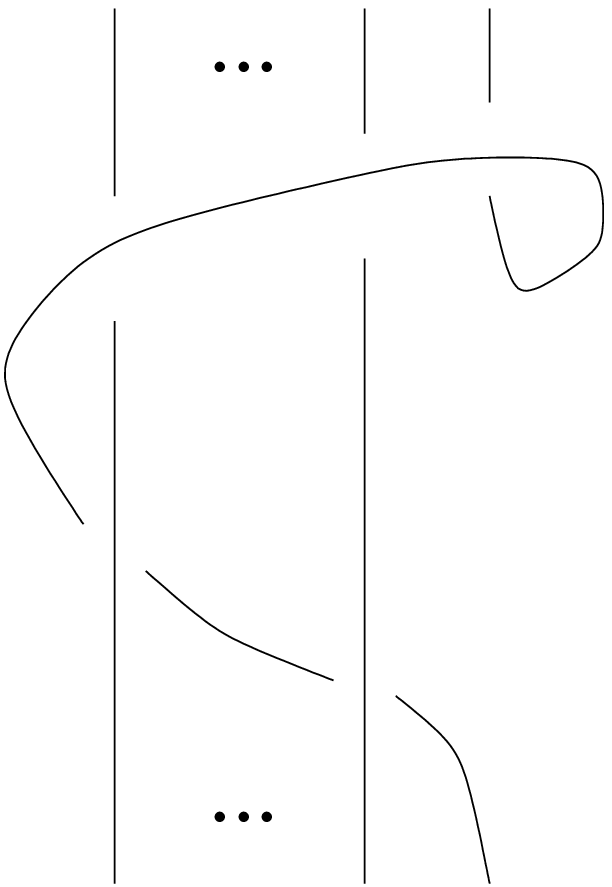}$$
For $i>0$ the equation in $K_{n+1}$,
$$h_n(\tau_n^i\otimes \oneun)h_n=Z_n^{(i)}\otimes h_1\ ,$$
defines a central element $Z_n^{(i)}$ in $K_{n-1}$.
We  consider the formal power series in $u^{-1}$,
$$Z_n(u)=\sum_{i\geq0}  Z_n^{(i)}  u^{-i}\ .$$
We have that
$$Z_n(u)\otimes h_1=h_n\left(\frac{u}{u-\tau_n}\otimes \oneun\right)h_n\ .$$
We denote by $Z_n(\mu,u)$ the series given by the action of
$Z_n(u)$ on the simple component of $K_{n-1}$ indexed by
$\mu$.

A canonical basis of $K_{n}$ is given in Theorem \ref{kbase}.
 Let $\Lambda$ be an up and down tableau, whose length
is $n+1$, and such that
$$\Lambda_{n-1}=\Lambda_{n+1}=\m\text{ and }\Lambda_n=\l\ .$$
Write the products $h_nq_\Lambda$ and $q_\Lambda h_n$
in the canonical basis.
$$h_n q_\Lambda=\sum_{\Xi_{n+1}=\mu} h_n(\Xi,\Lambda)\mathfrak{a}_\Xi
\mathfrak{b}_\Lambda\;\;\;\text{ and }\;\;\;
q_\Lambda h_n=\sum_{\Xi_{n+1}=\mu} h'_n(\Lambda,\Xi)\mathfrak{a}_\Lambda
\mathfrak{b}_\Xi\ .$$
Let $J_\Lambda$ be the set of up and down tableaus 
$\Xi=(\Xi_1,\dots,\Xi_{n+1})$ such that
$\Xi_m=\Lambda_m$ for every $m\neq n$.
 By considering $q_\Xi h_n q_\Lambda$ 
(resp. $q_\Lambda h_n q_\Xi$), we get 
$$h_n q_\Lambda=\sum_{\Xi\in J_\Lambda} h_n(\Xi,\Lambda)\mathfrak{a}_\Xi
\mathfrak{b}_\Lambda\;\;\; \text{ and }\;\;\;
q_\Lambda h_n=\sum_{\Xi\in J_\Lambda} h'_n(\Lambda,\Xi)\mathfrak{a}_\Lambda
\mathfrak{b}_\Xi\ .$$
Using the  three lemmas below the proof can be accomplished as follows.
From Lemma \ref{lemm1}, Lemma \ref{lemm2} and (\ref{serq1}) 
 we have 
 $$\frac{\la\l\ra}{\la\m\ra}=\mathrm{res}_{u=b}\frac{Z(\mu,u)}{u}=
\mathrm{res}_{u=b}\frac{Q(\mu,u)}{u}\ $$
The required formula follows now from (\ref{recqb}).
\end{proof}

\begin{lem}\label{lemm1}
 One has $$
h_n(\Lambda,\Lambda)=h'_n(\Lambda,\Lambda)=\frac{\la \l \ra}{\la \m \ra}\ .$$
\end{lem}
\begin{proof}
Let $\Lambda'$ be as usual obtained by removing the last term in the sequence 
$\Lambda$.
We have 
$$q_\Lambda=
\mathfrak{a}_{\Lambda}
\mathfrak{b}_{\Lambda}
=\frac{\la \m \ra}{\la \l \ra}\,
(\mathfrak{a}_{\Lambda'}\otimes 1)
(\tilde y_\mu\otimes h_{\l/\mu})(\mathfrak{b}_{\Lambda'}\otimes 1)$$
We will obtain the diagonal term $h_n(\Lambda,\Lambda)$ by computing the 
quantum trace of
$h_nq_\Lambda$.
$$
h_n(\Lambda,\Lambda)\langle \mu\rangle=\langle h_nq_\Lambda\rangle
=\frac{\langle \mu\rangle}{\langle \l\rangle}\langle T \rangle\ ,
$$
where $T$ is represented in the picture below.\\[-8pt]
$$T=\ \ \psdiag{10}{30}{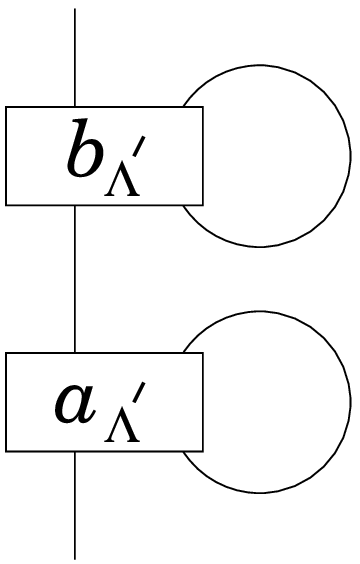}$$
Clearly, $\la T\ra=\la\l\ra^2\la\mu\ra^{-1}$ and we get the required 
result for
$h_n(\Lambda,\Lambda)$; $h'_n(\Lambda,\Lambda)$ is
calculated  similarly.
\end{proof}
Let us denote by
  $\mathrm{res}_{u=b}\ \frac{Z(\mu,u)}{u}$  the residue
of $\frac{Z(\mu,u)}{u}$ at $u=b$, i.e. the coefficient by $(u-b)^{-1}$ 
in the Laurent expansion of this function
 in the neighborhood of the point $u=b$.

\begin{lem}\label{lemm2}
$$h_n(\Lambda,\Lambda)=\mathrm{res}_{u=b}\ \frac{Z(\mu,u)}{u}$$
\end{lem}
\begin{proof}
We have that 
$$h_n(Z_n(u)\otimes \one\!{}_2) u^{-1}=h_n((u-\tau_n)^{-1}\otimes \oneun )h_n\ .$$
 Multiplying on the left by $q_\Lambda$, we can express the above 
formula  in the canonical basis.
 We denote by
$\zeta_\Lambda$ the diagonal term of index $\Lambda$.
 From the left hand side, we get 
$$\zeta_\Lambda=h_n(\Lambda,\Lambda)Z(\mu,u)u^{-1}\ .$$ 
Let us compute the coefficient from  the right hand side.
$$\begin{array}{rcl}
\zeta_\Lambda q_\Lambda&=&
q_\Lambda h_n((u-\tau_n)^{-1}\otimes \oneun )h_nq_\Lambda
\\
&=&\displaystyle \sum_{L\in J_\Lambda}\sum_{\Xi\in J_\Lambda}
h'_n(\Lambda,L)h_n(\Xi,\Lambda)\mathfrak{a}_\Lambda \mathfrak{b}_L
((u-\tau_n)^{-1}\otimes \oneun )
\mathfrak{a}_\Xi \mathfrak{b}_\Lambda\end{array}\ .$$
 In the the above sum, the term of indices
$L$ and $\Xi$ is zero unless $L=\Xi$, and in this case
 the action of $\tau_n$ multiplies
by the coefficient among $b_1,\dots,b_{2l+1}$ corresponding to
the Young diagram $\Xi_n$. These coefficients are distinct,
and we know that $h_n(\Lambda,\Lambda)=\frac{\langle\l\rangle}
{\langle\mu\rangle}$ is not zero. This implies
 that $Z(\mu,u)u^{-1}$ is a rational function in $u$,
whose residue at $u=b$ is equal to $h_n(\Lambda,\Lambda)$.
\end{proof}

\noindent The problem is now to compute the series $Z(\mu,u)$.
 We set 
\be\label{serq}Q_n(u)=Z_n(u)+\frac{\alpha^{-1}}{s-s^{-1}}-\frac{u^2}{u^2-1}\ ,
\ee
\be\label{serq1}
Q(\mu,u)=Z(\mu,u)+\frac{\alpha^{-1}}{s-s^{-1}}-\frac{u^2}{u^2-1}\ .\ee
\begin{lem}\label{lemQ}
\be \label{recqa}Q_{n+1}(u)\left(\frac{1}{u}-(s-s^{-1})^2\frac{\tau_n}{(u-\tau_n)^2}\right)
=
(Q_n(u)\otimes \oneun)
\left(\frac{1}{u}-(s-s^{-1})^2\frac{\tau_n^{-1}}{(u-\tau_n^{-1})^2}\right)
\ee
\be \label{recqb}Q(\mu,u)=\left(\frac{\alpha}{s-s^{-1}}+\frac{u\alpha}{u^2-1}\right)
\prod_j\frac{u-b_j^{-1}}{u-b_j}
\ee
\end{lem}
\begin{proof}
In the following computations, we will drop some $\oneun$.
This means that by drawing
the figures corresponding to these computation, one may have to add a vertical
string on the right in order to get  coherent equalities.
 For example, we write 
\be e_n^{-1}\tau_{n+1}=\tau_ne_n\ .\ee
From the above, using the skein relation, we obtain
\begin{eqnarray}
e_n^{-1}(u-\tau_{n+1})=(u-\tau_n)e_n^{-1}-(s-s^{-1})\tau_n(\one\!{}_{n+1}-h_n)\ ,\\
(u-\tau_n)e_n=e_n(u-\tau_{n+1})+(s-s^{-1})(\one\!{}_{n+1}-h_n)\tau_{n+1}\ .
\end{eqnarray}
This implies the following equalities for the formal series
\begin{eqnarray}
\label{un}\frac{1}{u-\tau_n}e_n^{-1}=e_n^{-1}\frac{1}{u-\tau_{n+1}}-(s-s^{-1})
\frac{\tau_n}{u-\tau_n}(\one\!{}_{n+1}-h_n)\frac{1}{u-\tau_{n+1}}\ ,
\end{eqnarray}\begin{eqnarray}
\label{deux}e_n\frac{1}{u-\tau_{n+1}}=\frac{1}{u-\tau_n}e_n+(s-s^{-1})
\frac{1}{u-\tau_n}(\one\!{}_{n+1}-h_n)\frac{\tau_{n+1}}{u-\tau_{n+1}}\ .
\end{eqnarray}
We also have the formula symmetric to (\ref{deux})
\begin{eqnarray}
\label{trois}\frac{1}{u-\tau_{n+1}}e_n=e_n\frac{1}{u-\tau_n}+(s-s^{-1})
\frac{\tau_{n+1}}{u-\tau_{n+1}}(\one\!{}_{n+1}-h_n)\frac{1}{u-\tau_n}\ .
\end{eqnarray}
We note that $\tau_n$ and $\tau_{n+1}$ commute, and that
$$h_n\tau_{n+1}^i =h_n\tau_n^{-i}\ .$$
Multiplying  (\ref{un}) on the left by $e_n$, we get
$$\begin{array}{lcl}
e_n\frac{1}{u-\tau_n}e_n^{-1}&=&
\frac{1}{u-\tau_{n+1}}-(s-s^{-1})
e_n\frac{1}{u-\tau_{n+1}}\frac{\tau_n}{u-\tau_n}\\
&&+(s-s^{-1})e_n\frac{u}{u-\tau_n}h_n\frac{1}{u-\tau_n^{-1}}
-(s-s^{-1})e_nh_n\frac{1}{u-\tau_n^{-1}}\ .
\end{array}$$
Using (\ref{deux}), (\ref{trois}) and the skein relations, we get
$$\begin{array}{lcl}
e_n\frac{1}{u-\tau_n}e_n^{-1}&=&
\frac{1}{u-\tau_{n+1}}-(s-s^{-1})\alpha^{-1}h_n\frac{1}{u-\tau_n^{-1}}
\\&&-(s-s^{-1})\frac{1}{u-\tau_n}e_n\frac{\tau_n}{u-\tau_n}
-(s-s^{-1})^2\frac{1}{u-\tau_n}\frac{\tau_{n+1}}{u-\tau_{n+1}}
\frac{\tau_n}{u-\tau_n}
\\&&+(s-s^{-1})^2\frac{1}{u-\tau_n}h_n\frac{\tau_n^{-1}}{u-\tau_n^{-1}}
\frac{\tau_n}{u-\tau_n}
+(s-s^{-1})\frac{u}{u-\tau_{n+1}}\alpha^{-1}h_n\frac{1}{u-\tau_n^{-1}}
\\&&-(s-s^{-1})^2u\frac{\tau_{n+1}}{u-\tau_{n+1}}\frac{1}{u-\tau_n}h_n
\frac{1}{u-\tau_n^{-1}}
+(s-s^{-1})^2u\frac{\tau_n^{-1}}{u-\tau_n^{-1}}h_n
\frac{1}{u-\tau_n}h_n\frac{1}{u-\tau_n^{-1}}\ .
\end{array}$$
We multiply on each side by $h_{n+1}$, and use the  relations 
$$h_{n+1}h_nh_{n+1}=h_{n+1}\text{ , }h_{n+1}e_nh_{n+1}=\alpha h_{n+1}
\text{ and }\tau_nh_{n+1}=h_{n+1}\tau_n\ .$$
$$\begin{array}{lcl}
\frac{Z_n(u)}{u}h_{n+1}&=&
\left(\frac{Z_{n+1}(u)}{u}-(s-s^{-1})\alpha^{-1}\frac{1}{u-\tau_n^{-1}}
\right.\\&&-(s-s^{-1})\alpha\frac{\tau_n}{(u-\tau_n)^2}
-(s-s^{-1})^2\frac{\tau_n}{(u-\tau_n)^2}
(Z_{n+1}(u)-\frac{\alpha-\alpha^{-1}}{s-s^{-1}}-1)
\\&&+(s-s^{-1})^2
\frac{1}{(u-\tau_n)^2(u-\tau_n^{-1})}
+(s-s^{-1})\frac{u\alpha^{-1}}{(u-\tau_n^{-1})^2}
\\&&-\left.(s-s^{-1})^2u\frac{\tau_n^{-1}}{(u-\tau_n)(u-\tau_n^{-1})^2}
+(s-s^{-1})^2\frac{\tau_n^{-1}}{(u-\tau_n^{-1})^2}Z_n(u)
\right)h_{n+1}\ .
\end{array}$$
The recursive formula (\ref{recqa}) can be deduced. 
It can be written
$$Q_{n+1}(u)=Q_{n}(u)\left(\frac{(u-\tau_n)^2}{(u-\tau_n^{-1})^2}
\frac{(u-s^{-2}\tau_n^{-1})(u-s^{2}\tau_n^{-1})}
{(u-s^{-2}\tau_n)(u-s^{2}\tau_n)}\right)\ .$$
Hence we have that
$$Q(\lambda,u)=Q(\mu,u)\left(\frac{(u-b)^2}{(u-b^{-1})^2}
\frac{(u-s^{-2}b^{-1})(u-s^{2}b^{-1})}
{(u-s^{-2}b)(u-s^{2}b)}\right)\ .$$
Recall that $b=\alpha s^{2{cn}(\l/\m)}$ is
the eigenvalue of $\tau_n$ corresponding to $q_{\Lambda'}$.
 The formula (\ref{recqb}) is then established recursively.
Note that 
$$Z(1,u)=\frac{u}{u-\alpha}
\left(\frac{\alpha-\alpha^{-1}}{s-s^{-1}}+1\right)\ .$$
Whence we have the formula for $Q(1,u)$.
\end{proof}

Using Theorem \ref{recdim}
we can deduce Wenzl's dimension formula \cite[Theorem 5.5]{Wbcd}.
Here $\lambda^{\!\vee}$ denote the transposed
Young diagram, so that $\lambda^{\!\vee}_j$ is the length of the
$j$th column of $\l$.
Let $n\in {\mathbb N}$ and $d\in {\mathbb Z}$, we set 
$$[y+d]=\frac{\a
s^d-\a^{-1}s^{-d}}{s-s^{-1}}, \;\;\;\; [n]=\frac{s^n-s^{-n}}{s-s^{-1}}.$$
\begin{thm}[Wenzl's formula]
\be \label{wen}\la \lambda\ra= 
\prod_{(j,j)\in\l}\frac{[y+\lambda_j-\lambda^{\!\vee}_j] 
+[hl(j,j)]}{[hl(j,j)]} 
\underset{i\neq j}{\prod_{(i,j)\in\l}}\frac{[y+d_\l(i,j)]}{[hl(i,j)]}\ee 
\end{thm}
Here, $hl(i,j)$ denote the hook-length of the cell 
$(i,j)$, i.e. $hl(i,j)=\l_i+\l^{\!\vee}_j-i-j+1$, 
and $d_\l(i,j)$ is defined by 
$$d_\l(i,j)=\left\{\begin{array}{lcl} 
\lambda_i+\lambda_j-i-j+1&\text{ if }&i\leq j\\ 
-\l^{\!\vee}_i-\l^{\!\vee}_j+i+j-1&\text{ if }&i>j\ . 
\end{array}\right. $$ 
If we define $d'_\l(i,j)$ by
$$d'_\l(i,j)=\left\{\begin{array}{lcl} 
\lambda_i+\lambda_j-i-j+1&\text{ if }&i<j\\ 
-\l^{\!\vee}_i-\l^{\!\vee}_j+i+j-1&\text{ if }&i\geq j , 
\end{array}\right.  $$
then we can write Wenzl's formula  as follows.
\be \label{wenzltwo}\la \lambda\ra= 
{\prod_{(i,j)\in\l}}\frac{\alpha^\frac{1}{2}s^{\frac{1}{2}d_\l(i,j)}
-\alpha^{-\frac{1}{2}}s^{-\frac{1}{2}d_\l(i,j)}}
{s^{\frac{1}{2}hl(i,j)}-s^{-\frac{1}{2}hl(i,j)}}
{\prod_{(i,j)\in\l}}\frac{\alpha^\frac{1}{2}s^{\frac{1}{2}d'_\l(i,j)}
+\alpha^{-\frac{1}{2}}s^{-\frac{1}{2}d'_\l(i,j)}}
{s^{\frac{1}{2}hl(i,j)}+s^{-\frac{1}{2}hl(i,j)}}\ee  
\begin{note}
It would be nice to interpret the above formula 
by decomposing $\lambda$ as a tensor product of two objects
in some bigger category.
\end{note}
\begin{proof} 
We will prove the formula (\ref{wenzltwo}).
 We first write the recursive formula (\ref{rec}) in a more convenient form.
We denote by $(i,\l_i)=(\l_{i'}^{\!\vee},i')$ the unique cell
in the skew diagram $\l/\m$.
\be \label{recvar}\begin{array}{rcl}
\frac{\la\l\ra}{\la\m\ra}&=&
\left(
\frac{\alpha s^{2\l_i-2i}-\alpha^{-1}s^{-2\l_i+2i}}
{s-s^{-1}}+1\right)\\
&&
\times \prod_{j<i}\frac{(\alpha s^{\l_i-i+\l_j-j+1}-
\alpha^{-1} s^{-\l_i+i-\l_j+j-1})
 (s^{\l_i-i-\l_j+j}- s^{-\l_i+i+\l_j-j})}
{(\alpha s^{\l_i-i+\l_j-j}-\alpha^{-1} s^{-\l_i+i-\l_j+j})
 (s^{\l_i-i-\l_j+j-1}- s^{-\l_i+i+\l_j-j+1})}\\
&&
\times \prod_{j'<i'}
\frac{(\alpha s^{-\l^{\!\vee}_{i'}+i'-\l^{\!\vee}_{j'}+j'-1}-
\alpha^{-1} s^{\l^{\!\vee}_{i'}-i'+\l^{\!\vee}_{j'}-j'+1})
 (s^{-\l^{\!\vee}_{i'}+i'+\l^{\!\vee}_{j'}-j'}- s^{\l^{\!\vee}_{i'}-i'-\l^{\!\vee}_{j'}+j'})}
{(\alpha s^{-\l^{\!\vee}_{i'}+i'-\l^{\!\vee}_{j'}+j'}-
\alpha^{-1} s^{\l^{\!\vee}_{i'}-i'+\l^{\!\vee}_{j'}-j'})
 (s^{-\l^{\!\vee}_{i'}+i'+\l^{\!\vee}_{j'}-j'+1}- s^{\l^{\!\vee}_{i'}-i'-\l^{\!\vee}_{j'}+j'-1})} 
\end{array}\ee
Here the first big product gives the contribution of the coefficients
$b_\xi$ corresponding to cells in the rows $1$ to $i-1$. Note that
some factors cancel if two among these rows have equal length.

We can write $\la\l\ra=\psi_\l(\alpha^{\frac{1}{2}},s^{\frac{1}{2}})
\psi'_\l (\alpha^{\frac{1}{2}},s^{\frac{1}{2}})$,
where $\psi_\l$ 
and $\psi'_\l$ satisfy the following recursive
formulas.
$$\begin{array}{rcl}
 \frac{\psi_\l(\beta,t)}{\psi_\m(\beta,t)}
&=& \frac{\beta t^{2\l_i-2i+1}-\beta^{-1} t^{-2\l_i+2i-1}}
{t-t^{-1}}
\\
&&
\times
\prod_{j<i}\frac{(\beta t^{\l_i-i+\l_j-j+1}-
\beta^{-1} t^{-\l_i+i-\l_j+j-1})
 (t^{\l_i-i-\l_j+j}- t^{-\l_i+i+\l_j-j})}
{(\beta t^{\l_i-i+\l_j-j}-\beta^{-1} t^{-\l_i+i-\l_j+j})
 (t^{\l_i-i-\l_j+j-1}- t^{-\l_i+i+\l_j-j+1})}\\
&&
\times \prod_{j'<i'}
\frac{(\beta t^{-\l^{\!\vee}_{i'}+i'-\l^{\!\vee}_{j'}+j'-1}-
\beta^{-1} t^{\l^{\!\vee}_{i'}-i'+\l^{\!\vee}_{j'}-j'+1})
 (t^{-\l^{\!\vee}_{i'}+i'+\l^{\!\vee}_{j'}-j'}- t^{\l^{\!\vee}_{i'}-i'-\l^{\!\vee}_{j'}+j'})}
{(\beta t^{-\l^{\!\vee}_{i'}+i'-\l^{\!\vee}_{j'}+j'}-
\beta^{-1} t^{\l^{\!\vee}_{i'}-i'+\l^{\!\vee}_{j'}-j'})
 (t^{-\l^{\!\vee}_{i'}+i'+\l^{\!\vee}_{j'}-j'+1}- t^{\l^{\!\vee}_{i'}-i'-\l^{\!\vee}_{j'}+j'-1})} 
\end{array}$$
$$\begin{array}{rcl}
 \frac{\psi'_\l(\beta,t)}{\psi'_\m(\beta,t)}
&=& \frac{\beta t^{-2\l^{\!\vee}_{i'}+2i'-1}+\beta^{-1} t^{2\l^{\!\vee}_{i'}-2i'+1}}
{t+t^{-1}}
\\
&&
\times
\prod_{j<i}\frac{(\beta t^{\l_i-i+\l_j-j+1}+
\beta^{-1} t^{-\l_i+i-\l_j+j-1})
 (t^{\l_i-i-\l_j+j}+t^{-\l_i+i+\l_j-j})}
{(\beta t^{\l_i-i+\l_j-j}+\beta^{-1} t^{-\l_i+i-\l_j+j})
 (t^{\l_i-i-\l_j+j-1}+t^{-\l_i+i+\l_j-j+1})}\\
&&
\times \prod_{j'< i'}
\frac{(\beta t^{-\l^{\!\vee}_{i'}+i'-\l^{\!\vee}_{j'}+j'-1}+
\beta^{-1} t^{\l^{\!\vee}_{i'}-i'+\l^{\!\vee}_{j'}-j'+1})
 (t^{-\l^{\!\vee}_{i'}+i'+\l^{\!\vee}_{j'}-j'}+t^{\l^{\!\vee}_{i'}-i'-\l^{\!\vee}_{j'}+j'})}
{(\beta t^{-\l^{\!\vee}_{i'}+i'-\l^{\!\vee}_{j'}+j'}+
\beta^{-1} t^{\l^{\!\vee}_{i'}-i'+\l^{\!\vee}_{j'}-j'})
 (t^{-\l^{\!\vee}_{i'}+i'+\l^{\!\vee}_{j'}-j'+1}+t^{\l^{\!\vee}_{i'}-i'-\l^{\!\vee}_{j'}+j'-1})} 
\end{array}$$
By induction, we can obtain the general formulas for
$\psi_\l(\beta,t)$ and $\psi'_\l(\beta,t)$.
$$\psi_\l(\beta,t)=
{\prod_{(i,j)\in\l}}\frac{\beta t^{d_\l(i,j)}
-\beta^{-1} t^{-d_\l(i,j)}}
{t^{hl(i,j)}-t^{-hl(i,j)}}$$
$$\psi'_\l(\beta,t)=
{\prod_{(i,j)\in\l}}\frac{\beta t^{d'_\l(i,j)}
+\beta^{-1} t^{-d'_\l(i,j)}}
{t^{hl(i,j)}+t^{-hl(i,j)}}$$
Whence we get  (\ref{wenzltwo}).
\end{proof}

The following proposition gives 
the quantum dimension formulas for the specializations corresponding
to the quantum groups of  B,C,D series (compare \cite{E}). 
 We consider here only partitions with at most $n$ rows.

\begin{pro}\label{Erl}
a) For  $\alpha=s^{2n}$ ($B_n$ specialization),
one has for a partition  $\lambda=(\l_1,...,\l_n)$
(which may have zero coefficients), 
$$\la\l\ra=\prod_{j=1}^n\frac{[n+\lambda_j-j+1/2]}{[n-j+1/2]}
\prod_{1\leq i< j\leq n}\frac{[2n+\lambda_i-i+\lambda_j-j+1]
[\lambda_i-i-\lambda_j+j]}
{[2n-i-j+1]
[j-i]}
$$
b) For $\alpha=s^{2n-1}$  ($D_n$ specialization),
one has for a partition  $\lambda=(\l_1,...,\l_n)$, 
$$\la\l\ra=
\prod_{1\leq i<j\leq n}\frac{[2n+\lambda_i-i+\lambda_j-j]
[\lambda_i-i-\lambda_j+j]}
{[2n-i-j]
[j-i]} \text{ if $\lambda_n=0$;}
$$
$$\la\l\ra=
2\prod_{1\leq i<j\leq n}\frac{[2n+\lambda_i-i+\lambda_j-j]
[\lambda_i-i-\lambda_j+j]}
{[2n-i-j]
[j-i]}\text{ if $\lambda_n\neq0$;}
$$
c) For $\alpha=-s^{2n+1}$  ($C_n$ specialization),
one has for a partition  $\lambda=(\l_1,...,\l_n)$, 
$$\la\l\ra=(-1)^{|\l|}
\prod_{j=1}^n\frac{[2n+2+2\lambda_j-2j]}{[2n+2-2j]}
\prod_{1\leq i<j\leq n}\frac{[2n+2+\lambda_i-i+\lambda_j-j]
[\lambda_i-i-\lambda_j+j]}
{[2n+2-i-j]
[j-i]}
$$
\end{pro}
\begin{note}
Observing that
$\langle \l\rangle_{\alpha,s}=\langle \l^\vee\rangle_{\alpha,-s^{-1}}=
\langle \l^\vee\rangle_{-\alpha^{-1},s}\ $,
we get formulas for the specializations which are {\em symmetric} to 
the above ones
(compare with \cite{E}).
For example, for  $\alpha=s^{-2n-1}$, 
one has  for a partition $\lambda$ whose first part $\l_1$ is  at most $n$,
and whose transposed partition is  
$\lambda^{\!\vee}=(\l_1^{\!\vee},...,\l_n^{\!\vee})$,
$$\la\l\ra=(-1)^{|\l|}
\prod_{j=1}^n\frac{[2n+2+2\lambda_j^{\!\vee}-2j]}{[2n+2-2j]}
\prod_{1\leq i< j\leq n}\frac{[2n+2+\lambda_i^{\!\vee}-i+\lambda_j^{\!\vee}-j]
[\lambda_i^{\!\vee}-i-\lambda_j^{\!\vee}+j]}
{[2n+2-i-j]
[j-i]}$$
\end{note}
\begin{proof}
Suppose that
$\l=(\l_1,\dots,\l_n)$,
and that $\mu$ is obtained from the Young diagram
$\l$ by removing one cell from the $i$th row,
 then we can write the recursive formula
 (\ref{recvar}) as follows.

\noindent If $i=n$ and $\l_n=1$, then
\be \label{recvar1}\begin{array}{rcl}
\frac{\la\l\ra}{\la\m\ra}&=&
\left(
\frac{\alpha s^{2-2n}-\alpha^{-1}s^{-2+2n}}
{s-s^{-1}}+1\right)\\
&&
\times \prod_{j<n}\frac{(\alpha s^{\l_n-n+\l_j-j+1}-
\alpha^{-1} s^{-\l_n+n-\l_j+j-1})
 (s^{\l_n-n-\l_j+j}- s^{-\l_n+n+\l_j-j})}
{(\alpha s^{\l_n-n+\l_j-j}-\alpha^{-1} s^{-\l_n+n-\l_j+j})
 (s^{\l_n-n-\l_j+j-1}- s^{-\l_n+n+\l_j-j+1})}\ ,
\end{array} \ee
else
\be \label{recvar2}\begin{array}{rcl}
\frac{\la\l\ra}{\la\m\ra}&=&
\left(
\frac{\alpha s^{2\l_i-2i}-\alpha^{-1}s^{-2\l_i+2i}}
{s-s^{-1}}+1\right)\\
&&
\times \prod_{\underset{j\neq i}{1\leq j\leq n}}
\frac{(\alpha s^{\l_i-i+\l_j-j+1}-
\alpha^{-1} s^{-\l_i+i-\l_j+j-1})
 (s^{\l_i-i-\l_j+j}- s^{-\l_i+i+\l_j-j})}
{(\alpha s^{\l_i-i+\l_j-j}-\alpha^{-1} s^{-\l_i+i-\l_j+j})
 (s^{\l_i-i-\l_j+j-1}- s^{-\l_i+i+\l_j-j+1})}\\
&&
\times \frac{
(s-s^{-1})}
{
(\alpha s^{2\l_i-2i-1}-\alpha^{-1} s^{-2\l_i+2i+1})
}
\times \frac{(\alpha s^{-n+\l_i-i}-\alpha^{-1} s^{n-\l_i+i})
}
{
( s^{n+\l_i-i}- s^{-n-\l_i+i})}\end{array}\ee

\noindent In the a) case, $\alpha=s^{2n}$, 
if $i=n$ and $\l_n=1$, we get
$$
 \frac{\la\l\ra}{\la\m\ra}=
(s+1+s^{-1})
\prod_{1\leq i<  n}\frac{[n+\lambda_j-j+2]
[1-n-\lambda_j+j]}
{[n+\lambda_j-j+1]
[-n-\lambda_j+j]}\ ,
$$
else
$$
\frac{\la\l\ra}{\la\m\ra}=\frac{[n+\lambda_i-i+1/2]}
{[n+\lambda_i-i-1/2]}
\prod_{\underset{j\neq i}{1\leq j\leq n}}
\frac{[2n+\lambda_i-i+\lambda_j-j+1]
[\lambda_i-i-\lambda_j+j]}
{ [2n+\lambda_i-i+\lambda_j-j]
[\lambda_i-i-\lambda_j+j-1]}\ .$$
The announced result follows.
\noindent In the b) case, $\alpha=s^{2n-1}$, 
if $i=n$ and $\l_n=1$, we get
$$
 \frac{\la\l\ra}{\la\m\ra}=
2
\prod_{1\leq i<  n}\frac{[n+\lambda_j-j+1]
[1-n-\lambda_j+j]}
{[n+\lambda_j-j]
[-n-\lambda_j+j]}\ ,
$$
else
$$
\frac{\la\l\ra}{\la\m\ra}=
\prod_{\underset{j\neq i}{1\leq j\leq n}}
\frac{[2n+\lambda_i-i+\lambda_j-j]
[\lambda_i-i-\lambda_j+j]}
{ [2n+\lambda_i-i+\lambda_j-j-1]
[\lambda_i-i-\lambda_j+j-1]}\ .$$
The formulas in b)  follow.
In the c) case, $\alpha=-s^{2n+1}$, 
if $i=n$ and $\l_n=1$, we get
$$
 \frac{\la\l\ra}{\la\m\ra}=
-(s^2+s^{-2})
\prod_{1\leq i<  n}\frac{[n+3+\lambda_j-j]
[1-n-\lambda_j+j]}
{[n+2+\lambda_j-j]
[-n-\lambda_j+j]}\ ,
$$
else
$$
\frac{\la\l\ra}{\la\m\ra}=-
\frac{[2n+2+2\lambda_i-2i]}
{[2n+2\lambda_i-2i]}
\prod_{\underset{j\neq i}{1\leq j\leq n}}
\frac{[2n+2+\lambda_i-i+\lambda_j-j]
[\lambda_i-i-\lambda_j+j]}
{ [2n+1+\lambda_i-i+\lambda_j-j]
[\lambda_i-i-\lambda_j+j-1]}\ .$$
The formula follows.
\end{proof}

\section{Formulas for the idempotents and the non generic case}
We conclude this paper with a summary  of  conditions needed
to define our minimal idempotents $\tilde p_t$ and 
$\tilde y_\l$. 
 This is of importance  in  the 
non generic case where $\a$ and $s$ are roots of unity 
and some of quantum integers $[m]$ are non invertible. 

 Let $\mu$ is a Young diagram
obtained from $\l$ by removing one cell.
From Corollary \ref{branch2} we get a formula for
the idempotent $\tilde y_\l$. 
 Here we denote by $\hat y_\lambda\in K_{\Box_\l}$ any lifting
of $y_\lambda\in H_{\Box_\l}$.
 We have
$$\hat y_\lambda(\tilde y_\m\otimes \one\!{}_1)\hat y_\lambda=
\sum_{\genfrac{}{}{0pt}{2}{\mu\subset \nu}{|\nu|=|\l|}}
\hat y_\lambda(\tilde y_\m\otimes \one\!{}_1)\tilde y_\nu
(\tilde y_\m\otimes \one\!{}_1)\hat y_\lambda
+\sum_{\genfrac{}{}{0pt}{2}{\nu\subset \mu}{|\nu|=|\l|-2}}
\hat y_\lambda\tilde y_{(\mu,\nu)}\hat y_\lambda\ .$$
Note that the first sum on the right hand side
lives in the Hecke summand of $K_{\Box_\l}$,
hence in each term we can replace $\hat y_\l$
by $\tilde y_\l$. By using Lemmas \ref{orth} and \ref{abs} we get: 
\be \tilde y_\lambda=
\hat y_\lambda(\tilde y_\m\otimes \one\!{}_1)\hat y_\lambda-
\sum_{\genfrac{}{}{0pt}{2}{\nu\subset \mu}{|\nu|=|\l|-2}}
\hat y_\lambda\tilde y_{(\mu,\nu)}\hat y_\lambda\ee
From the above formula we obtain a minimal idempotent $\tilde y_\l$
if the following three conditions are satisfied: 
\begin{itemize}
\item the quantum integer $[m]$ is non zero  for
any $m<\l_1+\l^\vee_1$;
\item the idempotent $\tilde y_\mu$ is defined for {\it some} $\mu\subset \l$,
$|\mu|=|\l|-1$;
\item
the coefficient $\frac{\la\mu\ra}{\la\nu\ra}$ 
given in Theorem \ref{recdim} is  nonzero for any  $\nu
\subset\mu\subset\l$, $|\nu|=|\mu|-1=|\l|-2$.
\end{itemize}

Let $t$ be a standard tableau, with shapes $\l(t)=\l$,
and $\l(t')=\mu$. Then from the general formula \ref{defsn}
for the section $s_n$ ($n=|\l|$), we get
\be \tilde p_t= 
\tilde p^+_{t'}  \;\hat p_t\; \tilde p^+_{t'}\ .\ee
 Here $\hat p_t$ is any lifting in $K_n$ of the path
idempotent $p_t\in H_n$. We obtain a minimal idempotent $\tilde p_t$
if the following three conditions are satisfied: 
\begin{itemize}
\item the quantum integer $[m]$ is non zero  for
any $m<\l_1+\l^\vee_1$;
\item the idempotent $\tilde y_\mu$ is defined;
\item
the coefficient $\frac{\la\mu\ra}{\la\nu\ra}$ 
given in Theorem \ref{recdim} is  nonzero for any  $\nu
\subset\mu\subset\l$, $|\nu|=|\mu|-1=|\l|-2$.
\end{itemize}

\bibliographystyle{amsplain} 

\begin{thebibliography}{10} 
\bibitem{AM} A.K.  Aiston, H.R. Morton, {\em Idempotents of Hecke algebras  
of type $A$}, Journal of Knot Theory and Ram. Vol. 7, 
No. 4 (1998), 463--487.
\bibitem{Aiston2} A. K.~Aiston, {\em A skein theoretic proof of the hook formula
for quantum dimension}, (preprint 1997). 
\bibitem{BB} A. Beliakova, C. Blanchet, {\em Modular categories of types B,
C and D}, Preprint math.QA/0006227 (2000)
\bibitem{bw} J.S. Birman, H. Wenzl, {\em Braids, 
link polynomials and a new algebra}, Trans. of AMS, Vol. 313, No. 1 (1989), 
249--273. 
\bibitem{Bhec} C. Blanchet, {\em Hecke algebras, modular 
categories and $3$-manifolds quantum invariants}, 
Topology, Vol. 39, No. 1 (2000), 193--223. 
\bibitem{BHMV} C. {Blanchet}, N. {Habegger}, G. {Masbaum} and 
P. {Vogel}, {\em Topological Quantum Field Theories derived from the 
Kauffman bracket}\/, Topology. Vol. 34, No. 4 (1995), 883--927. 
\bibitem{BW} J. Birman, H. Wenzl, {\em Braids, link polynomials
and a new algebra}, Trans. AMS, Vol. 313, No. 1 (1989), 249--273.
\bibitem{Br} R. Brauer, {\em On algebras which are connected with
the semisimple continuous groups}, Ann. of Math. Vol. 38 (1937), 854--872.
\bibitem{E} J. Erlijman, {\em New subfactors from braid group representations},
Trans. AMS. Vol. 350, No. 1 (1998) 185-211.
\bibitem{KRT} C. Kassel, M. Rosso, V. Turaev, {\em Quantum groups 
and knots invariants}\/, Panoramas et Synth\`eses No. 5, SMF (1997). 
\bibitem{Macdo} I. G. Macdonald, {\em Symmetric functions
and Hall polynomials}, Oxford Math. Monographs (1979).
\bibitem{MT} H. R. Morton and P. Traczyk, {\em Knots
and algebras}, Contribuciones Matematicas en homaje al professor
D. Antonio Plans Sanz de Bremond, E. Martin-Peinador and A. Rodez editors,
University of Saragoza (1990), 201--220.
\bibitem{Mu} J. Murakami, {\em The representations of the $q$-analogue
of Brauer's centralizer algebras and the Kauffman polynomial of links},
Publ. Res. Inst. Math. Sci. Vol. 26, No. 6  (1990), 935--945.
\bibitem{Naz} M. Nazarov, {\em Young orthogonal form for
Brauer's centralizer algebra}\/, J. of Algebra. Vol. 182 (1996), 664--693.
\bibitem{RW} A. Ram, H.~Wenzl, {\em Matrix units for centralizer
algebras}, J. of Algebra. Vol. 145 (1992), 378--395.
\bibitem{Tu} V. G. Turaev, {\em Quantum invariants 
of knots and $3$-manifolds}\/, De Gruyter Studies in Math.: 18 (1994).  
\bibitem{T1} V. G. Turaev, {\em Operator invariants of tangles,
and $R$-matrices}, Math. USSR Izv. Vol. 35,  No. 2 (1990), 411--443.
\bibitem{TW} V. Turaev, H. Wenzl, {\em Quantum invariants of 
$3$-manifolds associated with classical simple Lie algebras}, 
Int. J. of Math. Vol. 4, No. 2 (1993), 323--358. 
\bibitem{TW2} V. Turaev, H. Wenzl, {\em Semisimple and modular categories 
from link invariants}, Math. Ann. Vol. 309 (1997), 411--461. 
\bibitem{Wbcd} H.~Wenzl, {\em Quantum groups and subfactors of 
type B, C, and D}\/, Comm. Math. Phys. Vol. 133 (1990), 383--432. 
\bibitem{Wenzl} H.~Wenzl, {\em Hecke algebras of type $A_n$
and subfactors}\/, Invent.~Math. Vol. 92 (1988), 349--383.
\bibitem{WX} H.~Wenzl, {\em On the structure of Brauer's centralizer algebras},
Annals of Math. Vol. 128 (1988), 173--193.
\bibitem{Weyl} H. Weyl, {\em The Classical groups}, Princeton University Press
(1939).
\bibitem{Yo} Y.~Yokota, {\em Skeins and quantum $SU(N)$ 
invariants of $3$-manifolds}\/, Math. Ann. Vol. 307 (1997), 109--138. 
\end{thebibliography}
 
\end{document}